\documentclass[a4paper,12pt]{article}
\usepackage{amsmath,amsfonts,amssymb}
\newtheorem{conjecture}{Conjecture}
\begin{document}
\title{Une conjecture en th\'eorie des partitions}
\author{Michel Lassalle\thanks{Centre National de la Recherche Scientifique, France}\\
Ecole Polytechnique\\
91128 Palaiseau, France\\
e-mail: lassalle @ labri.u-bordeaux.fr}
 
\maketitle

\begin{abstract}
Nous pr\'esentons une conjecture qui se formule de mani\`ere \'el\'ementaire dans le cadre 
de la th\'eorie des partitions. 
\end{abstract}

\section{Introduction} \label{S:intro} 

Nous revenons dans cette Note sur une conjecture que nous avons \'et\'e conduit \`a formuler  
dans un pr\'ec\'edent article~\cite{La}. Cette conjecture s'exprime dans le cadre de la th\'eorie 
classique des partitions. Il est remarquable que sa formulation ne n\'ecessite que des notions 
extr\^emement \'el\'ementaires. Nous la pr\'esentons ici sous sa forme la plus naturelle,  
et nous explicitons plusieurs de ses cons\'equences.

\section{Notations}\label{S:not}

Une partition $\lambda$ est une suite d\'ecroissante finie d'entiers positifs. On dit que le nombre
 $n$ d'entiers non nuls est la longueur de $\lambda$. On note
$\lambda  = ( {\lambda }_{1},...,{\lambda }_{n})$ 
et $n = l(\lambda)$. On dit que 
$\left|{\lambda }\right| = \sum\limits_{i = 1}^{n} {\lambda }_{i}$
est le poids de $\lambda$, et pour tout entier $i\geq1$ que 
${m}_{i} (\lambda)  = card \{j: {\lambda }_{j}  = i\}$
est la multiplicit\'e de $i$ dans $\lambda$. On identifie $\lambda$ \`a son diagramme de Ferrers 
$\{ (i,j) : 1 \le i \ \le l(\lambda), 1 \ \le j \ \le {\lambda }_{i} \}$. 
On pose 
\[{z}_{\lambda }  = \prod\limits_{i \ge  1}^{} {i}^{{m}_{i}(\lambda)} {m}_{i}(\lambda) !  .\]
	 
Nous avons introduit dans~\cite{La} la g\'en\'eralisation suivante du coefficient binomial classique. 
Soient $\lambda$ une partition et $r$ un entier $\geq0$. On note 
$\genfrac{\langle}{\rangle}{0pt}{}{\lambda}{r}$ 
le nombre de fa\c cons dont on peut choisir r points 
dans le diagramme de $\lambda$ de telle sorte que 
\textit{au moins un point soit choisi sur chaque ligne de $\lambda$}.
	
Soient $X$ une ind\'etermin\'ee et $n$ un entier $\geq1$. On note d\'esormais 
\[{(X)}_{n }  =  X (X +1) ... (X + n-1)\] 
\[{[X]}_{n }  =  X (X - 1) ... (X -n +1)  .\]	
les coefficients hyperg\'eom\'etriques "ascendant" et "descendant" classiques. On pose 
\[\binom{X}{n}  =  {\frac{{[X]}_{n}}{n!}}  .\]
	
\section{Notre conjecture} \label{S:conj}
  
Il s'agit d'une g\'en\'eralisation de la propri\'et\'e classique suivante, qui est par exemple 
d\'emontr\'ee au Chapitre 1, Section 2, Exemple 1 du livre de Macdonald~\cite{Ma}.
Soit $X$ une ind\'etermin\'ee. Pour tout entier  $n\geq1$ on a 
\[\sum_{\left|{\mu }\right| = n} {(-1)}^{n-l(\mu)} {\frac{{X}^{l(\mu)}}{{z}_{\mu }}} =  \binom{X}{n},\]  
\[\sum_{\left|{\mu }\right| = n} {\frac{{X}^{l(\mu)}}{{z}_{\mu }}}  = \binom{X + n - 1}{n}.\]	
Ces deux relations sont \'equivalentes en changeant $X$ en $-X$. 
	
Dans~\cite{La} nous avons formul\'e la conjecture suivante qui g\'en\'eralise ce r\'esultat. Soit $X$ une ind\'etermin\'ee. 
Pour tous entiers $n,r,s\geq1$  nous avons conjectur\'e l'identit\'e
\begin{multline*}
\sum_{\left|{\mu }\right| = n} {(-1)}^{r - l(\mu )} {\frac{\genfrac{\langle}{\rangle}{0pt}{}{\mu}{r}}
{{z}_{\mu }}} {X}^{l(\mu ) - 1} \left({ \sum_{i = 1}^{l(\mu )} {({\mu }_{i})}_{s} }\right)  \\
= (s - 1) ! \binom{n + s - 1}{n - r} \sum_{k = 1}^{\min (r,s)}  \binom{X - s}{r - k} \binom{s}{k} .
\end{multline*} 
	 
Mais il est bien connu (~\cite{Lo}, p.13) que la suite de polyn\^omes
$\{ {[X]}_{n} , n\ \ge  0 \}$ 
est du Òtype binomialÓ, c'est-\`a-dire qu'elle satisfait l'identit\'e
\[{[X + Y]}_{n} = \sum_{k \ge  0} \binom{n}{k} {[X]}_{n - k}{[Y]}_{k} .\]
	
Nous sommes ainsi conduits \`a pr\'esenter notre conjecture sous la forme suivante qui est 
beaucoup plus naturelle.

\begin{conjecture} \label{T:1}
	Soit $X$ une ind\'etermin\'ee. Pour tous entiers $n,r,s\geq1$ on a
\begin{multline*}
\sum_{\left|{\mu }\right| = n} {(-1)}^{r - l(\mu )} {\frac{\genfrac{\langle}{\rangle}{0pt}{}{\mu}{r}}{{z}_{\mu }}} 
{X}^{l(\mu ) - 1} \left({ \sum_{i = 1}^{l(\mu )} {({\mu }_{i})}_{s} }\right) \\
= (s - 1) ! \binom{n + s - 1}{n - r} \left[{\binom{X}{r} - \binom{X - s}{r}}\right] .
\end{multline*}	
Ou de mani\`ere \'equivalente
\begin{multline*}
\sum_{\left|{\mu }\right| = n} {\frac{\genfrac{\langle}{\rangle}{0pt}{}{\mu}{r}}{{z}_{\mu }}} 
{X}^{l(\mu ) - 1} \left({ \sum_{i = 1}^{l(\mu )} {({\mu }_{i})}_{s} }\right) \\
= (s - 1) ! \binom{n + s - 1}{n - r} \left[{\binom{X + r + s - 1}{r} - \binom{X + r - 1}{r}}\right] .
\end{multline*}
\end{conjecture}
	
L'\'equivalence est obtenue en  changeant $X$  en $-X$. Comme on a
\[\genfrac{\langle}{\rangle}{0pt}{}{\mu}{r}  = 0   \qquad\textrm{si } r > \left|{\mu }\right| ,\]
la Conjecture~\ref{T:1} est triviale pour $r>n$. Comme on a
\[\genfrac{\langle}{\rangle}{0pt}{}{\mu}{r}  = 0  \qquad\textrm{si } r  < l(\mu ) ,\]
la sommation au membre de gauche de la Conjecture~\ref{T:1} est limit\'ee aux partitions 
$\mu$ telles que $l(\mu)\leq r$.
Chacun des membres est ainsi un polyn\^ome en $X$ de degr\'e  $r-1$.
	
Comme on a
\[\genfrac{\langle}{\rangle}{0pt}{}{\mu}{\left|{\mu }\right|}  = 1 ,\]
la Conjecture~\ref{T:1} prend la forme suivante pour $r=n$.

\begin{conjecture} \label{T:2}
	Soit $X$ une ind\'etermin\'ee. Pour tous entiers $n,s\geq1$ on a
\[\sum_{\left|{\mu }\right| = n} {(-1)}^{n - l(\mu )} {\frac{{X}^{l(\mu ) - 1}}{{z}_{\mu }}}  \left({ \sum_{i = 1}^{l(\mu )} 
{({\mu }_{i})}_{s} }\right)  =  (s - 1) !  \left[{\binom{X}{n}  -  \binom{X - s}{n}}\right]  .\]
	Ou de mani\`ere \'equivalente	
\[\sum_{\left|{\mu }\right| = n}  {\frac{{X}^{l(\mu ) - 1}}{{z}_{\mu }}}  \left({ \sum_{i = 1}^{l(\mu )} 
{({\mu }_{i})}_{s} }\right)  =  (s - 1) !  \left[{\binom{X + n + s - 1}{n} - \binom{X + n - 1}{n}}\right] .\]
\end{conjecture}
	
La Conjecture~\ref{T:2} est v\'erifi\'ee pour $s=1$ car on retrouve dans ce cas la propri\'et\'e classique \'enonc\'ee au 
commencement de cette section.

\section{Commentaires} \label{S:comm}

Nous avons v\'erifi\'e la Conjecture~\ref{T:1} dans les cas particuliers suivants:
\begin{itemize}
      \item[(i)] pour $n\leq7$ avec $r$ et $s$ arbitraires (au moyen d'un calcul explicite),
      \item[(ii)] pour $s=1$ avec $n$ et $r$ arbitraires (c'est le th\'eor\`eme 1 de~\cite{La}),
      \item[(iii)] pour $r=1,2,3$ avec $n$ et $s$ arbitraires (voir ci-dessous).
\end{itemize}
	
En g\'en\'eral la Conjecture~\ref{T:1} se d\'ecompose en $r-1$ conjectures obtenues en identifiant les coefficients de 
${X}^{k} ( 0 \ \le  k \ \le  r - 1 )$ 
dans chaque membre.
	
Pour $k=0$ la sommation au membre de gauche est limit\'ee \`a la partition $(n)$ de longueur $1$. 
En identifiant les termes constants dans chaque membre, on obtient	
\[{(-1)}^{r} {\frac{\binom{n}{r}}{n}}  {(n)}_{s}  =  (s - 1) ! \binom{n + s - 1}{n - r}  \binom{- s}{r} .\]
Cette identit\'e est tr\`es facilement v\'erifi\'ee.

D'autre part on a le d\'eveloppement
\[\left[{\binom{X}{r} - \binom{X - s}{r}}\right]  = 
{\frac{1}{r !}} \left({r s {X}^{r - 1} - {\frac{r (r - 1)}{2}} s (r + s - 1) {X}^{r - 2}  + \textrm{etc...}}\right).\]
	
Le coefficient de ${X}^{r - 1}$ 
au membre de gauche de la Conjecture~\ref{T:1} correspond \`a la sommation sur les partitions $\mu$ de longueur $r$. 
Comme  on a 
\[\genfrac{\langle}{\rangle}{0pt}{}{\mu}{l(\mu )}  = \prod_{i = 1}^{l(\mu )} {\mu }_{i} = \prod_{i \ge 1} {i}^{{m}_{i}(\mu )} ,\]	
en identifiant les coefficients de  ${X}^{r - 1}$ dans chaque membre, on obtient la

\begin{conjecture} \label{T:3}
	Pour tous entiers $n,r,s\geq1$  on a
\[(r - 1) !  \sum_{\begin{subarray}{1}\left|{\mu }\right| = n\\
l(\mu ) = r\end{subarray}} {\frac{\displaystyle\sum\limits_{i \ge  1}^{} {{m}_{i}(\mu ) (i)}_{s} }{\displaystyle\prod\limits_{i \ge  1}^{} {m}_{i}(\mu ) !}} 
=  s ! \binom{n + s - 1}{n - r} .\]

\end{conjecture}
	
Le coefficient de ${X}^{r - 2}$ au membre de gauche de la Conjecture~\ref{T:1} correspond \`a la sommation sur 
les partitions $\mu$ de longueur $r-1$.  On voit facilement qu'on a
\[\genfrac{\langle}{\rangle}{0pt}{}{\mu}{l(\mu ) + 1} =  {\frac{1}{2}} \left({\left|{\mu }\right| - l(\mu )}\right)
 \prod_{i = 1}^{l(\mu )} {\mu }_{i} .\]
	
En identifiant les coefficients de  ${X}^{r - 2}$ dans chaque membre, on retrouve la m\^eme Conjecture~\ref{T:3}, 
mais \'ecrite en remplacant $r$ par $r-1$.
	
La Conjecture~\ref{T:3} est v\'erifi\'ee pour $s = 0$ et $s = 1$, avec $n$ et $r$ arbitraires (voir
 \cite{La}, page 462). Pour $s = 0$ (resp $s = 1$) le membre de gauche est \'egal au nombre
  de ÒcompositionsÓ (c'est-\`a-dire de multi-entiers) de poids $n$ et de longueur $r$ 
(resp. ce nombre multipli\'e par $n/r$).

La Conjecture~\ref{T:3} est \'egalement v\'erifi\'ee pour $r = 1$ et $r = 2$, avec $n$ et $s$ arbitraires. 
Pour $r = 1$ elle devient l'identit\'e
\[{(n)}_{s}  = s ! \binom{n + s - 1}{s} .\]	
Pour $r=2$ elle s'\'ecrit
\[\sum_{i = 1}^{n - 1} {(i)}_{s} = s ! \binom{n + s - 1}{n - 2} .\]	
Cette identit\'e est v\'erifi\'ee car elle est \'equivalente \`a la propri\'et\'e classique suivante
\[\binom{N}{k}  = \sum_{i = 1}^{N - 1} \binom{i}{k - 1} .\]
		
Nous g\'en\'eralisons ce r\'esultat au moyen de la conjecture suivante.

\begin{conjecture} \label{T:4}
	Pour tous entiers $n,r,s\geq1$ on a.	
\[(r - 1) !  \sum_{\begin{subarray}{1}\left|{\mu }\right| = n\\
l(\mu ) = r\end{subarray}} {\frac{\displaystyle\sum\limits_{i \ge  1}^{} {{m}_{i}(\mu ) (i)}_{s} }{\displaystyle\prod\limits_{i \ge  1}^{} {m}_{i}(\mu ) !}}
=\sum_{i = 1}^{n - r + 1}  \binom{n - i - 1}{r - 2} {(i)}_{s}  .\]
\end{conjecture}
	
La Conjecture~\ref{T:4} permet de d\'emontrer facilement la Conjecture~\ref{T:3} au moyen d'une r\'ecurrence 
sur l'entier $r$. Nous l'avons v\'erifi\'ee explicitement pour  $r\geq n-7$ avec $n$ et $s$ 
arbitraires.


\begin{thebibliography}{9}
   \bibitem{La}
      M. Lassalle, \emph{Quelques conjectures combinatoires 
      relatives \`a la formule classique de Chu -Vandermonde}, Adv. in Appl. Math. 
      \textbf{21} (1998), 457--472.
   \bibitem{Lo}
      D. Loeb, G.C. Rota, \emph{Formal power series of logarithmic type}, 
      Adv. in Appl. Math. \textbf{75} (1989), 1--118.
   \bibitem{Ma}
      I.G. Macdonald, \emph{Symmetric functions and Hall polynomials}, Clarendon Press,
      second edition, Oxford, 1995.
\end{thebibliography}
\end{document}